\documentclass{amsart}
\usepackage{amssymb, mathbbol}
\usepackage{amsfonts, a4wide, amscd}
\usepackage{amsmath, amsthm, dsfont}
\usepackage[all]{xy}

\newtheorem{Theorem}{Theorem}[section]
\newtheorem{Lemma}[Theorem]{\sc Lemma}
\newtheorem{Definition}[Theorem]{\sc Definition}

\newtheorem{Remark}[Theorem]{\sc Remark}

\def\C{{\mathbb{C}}}
\def\L{{\mathbb{L}}}
\def\Z{{\mathbb{Z}}}
\def\A{{\mathbb{A}}}
\def\P{{\mathbb{P}}}

\def\Rcal{{\mathcal{R}}}

\def\m{{\mathcal{M}}}
\def\bz{{\overline{\mathcal{Z}}}}
\def\Om{{\Omega^*}}

\def\om{{\Omega_*}}
\def\o{{\Omega}}
\def\baro{{\overline{\Omega}}}
\def\pf{{\noindent{\sc Proof: }}}

\newcommand\ra{\rightarrow}
\newcommand{\bt}{\begin{Theorem}}
\newcommand{\et}{\end{Theorem}}
\newcommand{\bl}{\begin{Lemma}}
\newcommand{\el}{\end{Lemma}}
\newcommand{\bd}{\begin{displaymath}}
\newcommand{\ed}{\end{displaymath}}
\newcommand{\ba}{\begin{align*}}
\newcommand{\ea}{\end{align*}}

\begin{document}

\title{Algebraic Cobordism of Classifying Spaces}
\author{Dinesh Deshpande}
\address{Department of Pure Mathematics and Mathematical Statistics, University of Cambridge, Cambridge, CB3 0WB, UK}
\email{D.Deshpande@dpmms.cam.ac.uk}

\maketitle{}

\begin{abstract} We define algebraic cobordism of
classifying spaces, $\o^*(BG)$ and $G$-equivariant algebraic
cobordism $\o^*_G(-)$ for a linear algebraic group $G$. We prove
some properties of the  
coniveau filtration  on algebraic cobordism, denoted $F^j(\o^*(-))$,
which are required for the definition to work. We show that
$G$-equivariant cobordism satisfies the localization exact
sequence. We calculate 
$\o^*(BG)$ for algebraic groups over the complex numbers
corresponding to classical
Lie groups $GL(n)$, $SL(n)$, $Sp(n)$, $O(n)$ and $SO(2n+1)$. We also
calculate $\o^*(BG)$ when $G$ is a 
 finite abelian group. A finite non-abelian group for which we calculate
 $\o^*(BG)$ is the quaternion group of order 8. In all the above cases
 we check that $\o^*(BG)$ is isomorphic to $MU^*(BG)$.  \end{abstract}


\section{Introduction}

Complex cobordism $MU^*(-)$ is a cohomology theory represented
by the Thom spectrum $MU$ in the classical stable homotopy category.
Voevodsky has defined a version of algebraic cobordism, 
$MGL^{*,*}(-)$ for arbitrary schemes which can be found in
\cite{voe98}. We call $MGL^{*,*}(-) $ as motivic
cobordism. $MGL^{*,*}$ is a bigraded cohomology theory represented
by the motivic Thom spectrum $MGL$ in the Morel-Voevodsky motivic
homotopy category.  
Levine and Morel have constructed  algebraic cobordism
$\Om(-)$ \cite{levmor07}, based on the arguments of Quillen in
\cite{qui71}. For a quasi-projective schemes $X$, the ring $\Om(X)$ is
constructed 
as \begin{center} $\Om(X) = \{[f: Y \ra X]\}/(relations)$ \end{center}
where $f$ is a projective morphism of pure codimension from a smooth
variety $Y$ to $X$. Some 
details of this construction are given in section
2. Levine has recently showed that for smooth
quasi-projective schemes, the definition of algebraic
cobordism $\Om(-)$ agrees with the definition of motivic cobordism
$MGL^{2*,*}(-)$ \cite{lev08}.  

It turns out that the algebraic cobordism of a field is isomorphic to the
Lazard ring, which is isomorphic to the complex cobordism of a point
(\cite{ada95}, 2.8). This naturally leads to the question of finding
varieties whose algebraic cobordism is isomorphic to the complex cobordism of
 the realisation of the variety over the complex field. Classifying spaces
 of groups is a 
good class of varieties to consider. As a topological space,
classifying spaces are infinite dimensional for most groups. So one
needs to carefully define an analogue of classifying spaces in
algebraic geometry. Burt Totaro has defined and studied Chow rings of
classifying spaces of linear algebraic groups in \cite{tot99}. Let $G$ be a linear algebraic group. Let $G$ act on $\A^N$, i.e. consider an $N$-dimensional representation of $G$. Let $S$ be a (Zariski) closed subset of $\A^N$ of codimension $j$ such that $G$ acts freely on $\A^N-S$. Totaro defined Chow rings of the classifying space of $G$ by showing that $CH^i(\frac{\A^N-S}{G})$ is independent of the choices of $S$ and $\A^N$ whenever $i$ is less than $j$.
He has
also defined the cycle class map \bd CH^*(X) \ra MU^*(X) \otimes_{MU^*}
\Z \ed in \cite{tot97}. In all the examples in \cite{tot99}, it turns out
that the cycle class map $CH^*(BG) \ra MU^*(BG) \otimes_{MU^*} \Z$ is an
isomorphism. Algebraic cobordism relates to Chow rings in exactly
similar fashion \cite{levmor07} , i.e.\bd  CH^*(X) \cong \Om(X)
\otimes_{\Om} \Z . \ed
 Hence it is very natural to expect that the canonical map
from algebraic cobordism to complex cobordism is an isomorphism
whenever the cycle class map is an isomorphism. This paper is an attempt to
check \bd \Om(BG) \cong MU^*(BG)  \ed  for a few groups $G$ such as
classical Lie groups, products of 
finite cyclic groups and the quaternion group. We define $\Om(BG)$ for a linear algebraic groups $G$. Unlike Chow groups, $\Omega^i(\frac{\A^N-S}{G})$ are not independent of the choice of $S$. A reason being, $CH^i(X)$ is zero when $i < 0$ whereas $\Omega^i(X)$ is non-zero even when $i < 0$. Instead, we show that quotients of $\Omega^i(\frac{\A^N-S}{G})$ by the coniveau filtration are independent of the choices of $S$ and $\A^N$ [Theorem \ref{independent}]. Similar to
Edidin-Graham \cite{edigra98}, we define $G$-equivariant algebraic
cobordism for a $G$-scheme $X$. We show that $\o^*_G(-)$ satisfies the
localization exact sequence.  

Let $BP^*(-)$ be the Brown-Peterson cohomology for a prime $p$. Yagita
has defined $\Om_{BP}(-)$, the algebraic Brown-Peterson theory in the
section 8 of \cite{yag05} and has checked $\Om_{BP}(BG) \cong
BP^*(BG)$ for a few groups. We check that our result implies Yagita's result. $MGL^{*,*}(BGL(n))$ has been calculated in
\cite{hukri01}.

The structure of this chapter is as follows. Section 2 gives a brief
introduction to the theory $\Om(-)$. In section 3 we define algebraic
cobordism of classifying spaces similar to as done in
\cite{tot99}. This definition uses the coniveau filtration on cobordism
ring. We will also prove some properties of the coniveau filtration in this
section. Computations of $\Om(BG)$ for some reductive algebraic groups
$G$ over complex numbers corresponding to the classical Lie groups are
shown in  
section 4. Computations for products of the finite cyclic groups $G$ 
are done in 
section 5. Section 6 gives a way of calculating $\Om(X)$ when $X$ can
be partitioned into affine spaces. We define equivariant algebraic
cobordism in section 7 and calculate $\Om(BQ)$ for the group of
quaternions $Q$ in section 8. \newline 

{\it Acknowledgements}: I thank my PhD supervisor Burt Totaro for introducing me to this subject and for numerous interesting discussions.

\section{Brief Introduction to Algebraic Cobordism}

In this section, we will give two equivalent definitions of $\Om(-)$ extracted
from the chapter 2 of \cite{levmor07} and from
\cite{levpan09}and state some of the properties of $\Om(-)$.  

Throughout this chapter, assume that $k$ is a field of characteristic
zero and $k$ admits 
resolution of singularities. Let $Sch/k$ denote the category of
equidimensional quasi-projective schemes of finite type over $k$. By a
scheme, we always mean an element of $Sch/k$. Let
$Sm/k$ denote a full subcategory of $Sch/k$ whose objects are also
smooth. Let $X \in Sch/k$.  
 As said
in the introduction, \begin{center}$\om(X) = \{$bordism
  cycles$\}/(relations). $\end{center}  Define a bordism cycle on
$X$ to be a family $[f: Y \ra X, L_1, L_2, ... , L_r]$ where $Y$ is in
$Sm/k$, $f$ is projective morphism and $L_i$'s are line 
bundles over $Y$.  Define the degree of such a bordism cycle to be
$dim_k(Y) - r$.  Let $\mathcal{Z}_*(X)$ be the free abelian
group  generated by the isomorphism classes of the bordism cycles. 

We will now start imposing relations on $\mathcal{Z}_*(X)$ so that we
will have a working definition of Chern classes of vector bundles and
the formal group law to be the universal formal group law. Firstly we will
impose the 
dimensional constraint.  Let $\mathcal{R}_*^{dim}(X) \subset
\mathcal{Z}^*(X)$ be the subgroup generated by  cycles of the form $[Y
\ra X, L_1,..., L_r]$ where $dim_k(Y) < r$.  Define, 
\bd \bz _*(X) = \mathcal{Z}_*(X)/\mathcal{R}_*^{dim}(X). \ed

On $\bz _*(X)$ we will now impose the relation which will help us in
defining the Chern classes. Let $\mathcal{R}_*^{sect}(X) \subset \bz
_*(X)$ be the subgroup generated by all the elements of the form $([Y
\ra X, L] - [Z \hookrightarrow Y \ra X])$, where $L$ is a line bundle
over $Y$, $s: Y \ra L$ is a section transversal to the zero section and $Z
\hookrightarrow Y$ is a closed subvariety of the zeroes of $s$. Note that
$Z$ is smooth as $s$ is transversal to the zero section. Define, \bd
\m_*(X) = \bz_*(X)/\mathcal{R}_*^{sect}(X) . \ed 
 
Now we are ready to define the algebraic pre-bordism $\baro_*(X)$. We
say the two cycles $[f: Y \ra X]$ and $[g: Z \ra X]$ are {\it
  elementary bordant} if there exist $W \in Sm/k$ and a projective
morphism $h: W \ra X \times \P^1$ 
transversal to $X \times \{0\}$ and $X \times \{1\}$ such that $[h:
h^{-1}(X \times \{0\}) \ra X]$ is isomorphic to $[f: Y \ra X]$ and
$[h: h^{-1}(X \times \{1\}) \ra X]$ is isomorphic to $[g: Z \ra
X]$. Let $\Rcal_*^{cob}(X) \subset \m_*(X)$ be the subgroup generated by
elements of the form $([f: Y \ra X] - [g: Z \ra X])$ where $[f]$ and
$[g]$ are elementary bordant to each other. Define {\it Algebraic
  Pre-bordism} by, \bd \baro_*(X) = \m_*(X)/\Rcal_*^{cob}(X). \ed 

Given a line bundle $L$ over $X$, we define the Chern class homomorphism
$c_1(L) : \baro_*(X) \ra \baro_{*-1}(X) $ by $c_1(L)([f: Y \ra X]) =
[f: Y \ra X, f^*(L)]$. Note that $\Rcal_*^{dim}(X)$ forces $c_1(L)$ to
be locally nilpotent for all $L$. i.e. given any cycle $[f: Y \ra X]$,
there exist $n>o$ such that $(c_1(L))^n([f: Y \ra X]) = 0$.  

Now we will enforce the formal group law on $\baro_*(X)$, which will
give us the theory of algebraic bordism. Let \begin{center}$F(x,y) =
  x + y + \Sigma_{i,j>0} a_{ij}x^iy^j$\end{center} be the universal
formal group law with  coefficients in the Lazard ring $\L_*$. $\L_*$
is the quotient of the formal power series ring $\Z[[a_{11}, a_{12},...]]$
by the relations enforced via properties of $F(-,-)$ (\cite{ada95},
2.6). $\L_*$ is a graded ring and the degree of $a_{ij}$ is $(i+j-1)$. Let
$\Rcal_*^{FGL} \subset \L_* \otimes \baro_*(X)$ be the subgroup
generated by 
elements of the form $([F(c_1(L),c_1(M))(x)] - [c_1(L \otimes
M)(x)])$; where $L$, $M$ are line bundles over $X$ and $x \in
\baro_*(X)$. Local nilpotency of $c_1(-)$ ensures that $\Rcal_*^{FGL}$
is well defined. We define the {\it algebraic bordism} of $X$ by
\begin{align} \om(X) = 
\L_* \otimes \baro_*(X) / <\L_*.\Rcal_*^{FGL}> \end{align} 

The theory of algebraic bordism can also be constructed by allowing
the `double point degeneration' in the elementary bordism relations
\cite{levpan09}. A formulation is as follows.

Let $Y$ be in $Sch/k$. A morphism $\pi: Y \ra \P^1$ is a double point
degeneration over $0 \in \P^1$ if
\begin{align*}
  \pi^{-1}(0) = A \cup B
\end{align*}
where $A$ and $B$ are smooth Cartier divisors intersecting
transversely in $Y$. The intersection $D = A \cap B$ is called the
double point locus of $\pi$ over $0 \in \P^1$. Let $N_{A/D}$ and
$N_{B/D}$ denote the normal bundles of $D$ in $A$ and $B$
respectively. Since $\mathcal{O}_D(A+B)$ is trivial, we have
\begin{align*}
  N_{A/D} \otimes N_{B/D} \cong \mathcal{O}_D.
\end{align*}
Hence the two projective bundles $\P(\mathcal{O}_D \oplus N_{A/D}) \ra
D$ and $\P(\mathcal{O}_D \oplus N_{B/D}) \ra D$ are isomorphic. Let
$\P(\pi) \ra D$ denote the either of the two.

Let $p_1$ and $p_2$ denote the projections to the first and the second
factor of $X \times \P^1$ respectively. Let $W$ be a smooth variety
and let $h: W \ra X \times \P^1 $ be a 
projective morphism. Suppose the composition
\begin{align*}
  h_2 := p_2 \circ h : W \ra \P^1
\end{align*}
is a double point degeneration over $0$. Let $[A \ra X]$ , $[B \ra X]$
and $[\P(h_2) \ra X]$ be the elements in $\mathcal{Z}_*(X)$ obtained
from the fiber $h_2^{-1}(0)$. Without loss of generality, say $1$ is a
regular value of $p_2 \circ h$ and say
$[p_1 \circ h: h^{-1}(X \times \{1\}) \ra X]$ is isomorphic to $[f: Y
\ra X]$. Define 
an associated double point relation over $X$ by
\begin{align*}
  [f: Y \ra X] - [A \ra X] - [B \ra X] + [\P(h_2) \ra X].
\end{align*}
Let $\mathcal{R}_*^{dpr}(X)$ be the subgroup generated by all 
associated double point relations over $X$. Define the \emph{double
  point bordism theory} $\omega_*(-)$ by
\begin{align}
  \omega_*(X) := \mathcal{Z}_*(X)/\mathcal{R}_*^{dpr}(X). 
\end{align}

Then there is a canonical isomorphism
\begin{align}
  \omega_*(X) \cong \om(X).
\end{align}

$\om(-)$ satisfies all the standard properties of an algebraic homology theory
such as the existence of  
push forwards along proper morphisms and the existence of pull backs under
smooth morphism. To be precise, $\om(X)$ is an oriented Borel-Moore
functor of geometric type on the category $Sch/k$. Axioms
defining such a functor are given in (\cite{levmor07}, 2.1.11, 2.2.9).

If $X$ is in $Sm/k$ and of pure dimension $n$,
then define the {\it algebraic cobordism} of $X$ by, \begin{align} \Om(X) :=
\o_{n-*}(X). \end{align}

Having constructed $\Om(-)$, we will now
state a few important properties of $\Om(-)$.

\subsection {Localization Sequence} 

 Let $X$ is in $Sch/k$. For $Z$, a closed subscheme of codimension $r$
 in $X$, let $i: Z \ra X$ and $j: X-Z \ra X$ be inclusions. Then we
 have an exact sequence, (\cite{levmor07}, 3.2.7)
\begin{align*}\xymatrix{ {\om(Z)} \ar[r]^{i_*} & {\om(X)} \ar[r]^{j^*}
    & {\om(X-Z)} \ar[r] & 0} \end{align*}
 if both $X$ and $Z$ are smooth, we can write this sequence in the
 cohomological notation as
\begin{align*}\xymatrix{ {\Omega^i(Z)} \ar[r]^{i_*} &
    {\Omega^{i+r}(X)} \ar[r]^{j^*} & {\Omega^{i+r}(X-Z)} \ar[r] & 0}
\end{align*}

\subsection{Homotopy Invariance}

 Let $X \in Sch/k$. Let the dimension of $X$ be $d$ and let $p: E \ra X$ be a vector bundle over
 $X$. Then $p_*: \om(X) \ra \Omega_{*+d}(E)$ is an isomorphism, (\cite{levmor07},
 3.6.3).

\subsection{Projective Bundle Formula}

 Let $X \in Sch/k$. If $E$ is a rank $n+1$ bundle over $X$, then there
 is an isomorphism 
 $\Phi_{X,E} : \oplus_{j=0}^n \Omega_{*-n+j}(X) \ra
 \om(\mathbb{P}(E))$. Here $\Phi_{X,E}$ is described as $\oplus \xi^j
 \circ q^*$, where $q: \mathbb{P}(E) \ra X$ is the projection and $\xi$
 denotes the  $c_1(\mathcal{O}(1))$ operator for a canonical line bundle
 $\mathcal{O}(1)$ over $\mathbb{P}(E)$, (\cite{levmor07}, 3.5.2). 

This implies that $\Om(\C P^n) = \Om[x]/(x^{n+1})$ where $x =
c_1(\mathcal{O}(1))(id)$ as \newline $\o^{n+1}(\C P^n) = 0$. 

\subsection{Universality}

(\cite{levmor07}, 1.1.2) gives a definition of an oriented cohomology
theory over the category $Sm/k$. Algebraic cobordism
satisfy these axioms and in fact is the universal oriented cohomology
theory (\cite{levmor07}, 7.1.3). Hence given any other theory $h^*(-)$
satisfying these axioms, there is a unique map between functors
$\Om(-) \ra h^*(-)$, which induces homomorphism of the graded rings
$\Om(X) \ra h^*(X)$ for all $X$ in $Sm/k$. In particular there
is a canonical graded ring homomorphism $\Om(X) 
\ra CH^*(X)$, and in fact $\Om(k) \ra CH^*(k) $ induces an
isomoprhism (\cite{levmor07}, 4.5.1) \begin{align*}CH^*(X) \cong
  \Om(X) \otimes_{\Om} \Z. \end{align*}

$\Om$ is isomorphic to the Lazard ring $\L^*$ (with cohomological
grading) which is generated over 
$\Z$ by the monomials in $a_{11}, a_{12}, a_{21}, ...$ where the degree of
$a_{ij}$ is $-(i+j-1)$. The map \bd \Om \ra CH^* \cong \Z \ed is the one
sending each $a_{ij}$ to zero. 

Given a smooth variety $X$, then $MU^{2i}(X(\C))$ satisfies oriented
cohomology theory axioms, where $(X(\C))$ is the topological space of
complex points over $X$. Hence by the universality there is a natural map
$\Om(X) \ra MU^{2*}(X(\C))$. This map factors via another natural map,
$\Om(X) \ra MGL^{2*,*}(X)$, i.e. \begin{align*}\Om(X) \ra MGL^{2*,*}(X) \ra
MU^{2*}(X(\C)).\end{align*} 

The composition of the isomorphism $CH^*(X) \cong
\Om(X) \otimes_{\Om} \Z $ with the natural map $\Om(X) \ra
MU^{2*}(X(\C))$ gives the cycle class map \begin{align*}CH^*(X) \ra MU^*(X)
\otimes_{MU^*} \Z \end{align*} constructed in \cite{tot97}. 

In a similar way $\om(-)$ is the universal functor among oriented
Borel-Moore functors over the category of varieties and there exists
maps such as \bd \om(X) \ra CH_*(X) \text{ ; } \om(X) \ra MU^{BM}_{2*}(X) \ed 
between the corresponding homology theories.

\section{Coniveau Filtration}\label{def}

In this section we will prove some properties of the coniveau
filtration. In particular, we will show that the $i$th Chern class of a
vector bundle lies in the $i$th level of this filtration. We will
define algebraic cobordism ring for classifying spaces of linear
algebraic groups. The idea involved is to extend Totaro's definition for Chow
rings given in \cite{tot99} to the case of algebraic cobordism. A
difference in the algebraic cobordism case is that we need to pass
through the coniveau filtration on algebraic cobordism. Vishik
\cite{vis07} has showed 
that the coniveau filtration is multiplicative.  

 For $X \in Sm/k $, define the coniveau filtration on
 $\o^iX$ by 
\begin{align*} F^j(\o^iX) := \{&x \in \o^iX | x 
 \text{ restricts to zero in } \o^i(X-S)
 \text { for some closed subspace }\\ &  S \subset X 
 \text{ of the codimension at least }  j \}.
 \end{align*} 

 That is, $x \in F^j(\o^iX)$ if and only if, there is a closed $S
 \subset X$ of codimension at least $j$ such that $i^*(x) = 0$ for the
 inclusion $i: X-S \hookrightarrow X$. Hence we get the filtration \bd
 \o^iX = F^0(\o^iX) \supset F^1(\o^iX) \supset F^2(\o^iX) \supset
 \cdots \ed 

\bl \label{pullback} Let $f: X \ra Y$ be a proper map and let $r =
dim(Y) - dim(X)$, then $f_*(F^j(\o^i(X))) \subset
F^{j+r}(\o^{i+r}(Y))$  and if $g: X \ra Y$ is a smooth morphism, then
$g^*(F^j(\o^i(Y))) \subset F^j(\o^i(X))$. \el 

\pf Let $g: X \ra Y$ be a smooth morphism, and $\alpha \in
F^j(\o^i(Y))$. Let $\alpha|_{Y-S_Y}$ be zero for some closed $S_Y
\subset Y$ of codimension $d$ which is greater equal $j$. Let $S_X
\subset X$ be the pull-back of $S_Y \subset Y$ along $g$. Hence
$codim(S_X \subset X) \geqq codim(S_Y \subset Y) \geqq j$. Now
consider the following commutative diagram, 

\begin{align}
\xymatrix{ {\o^{i-d}(S_Y)} \ar[r]^{i^Y_*} \ar[d]^{g^*} & {\o^i(Y)}
  \ar[r]^{j^{Y^*}} \ar[d]^{g^*} & {\o^i(Y-S_Y)} \ar[r] & {0} \\ 
{\o^{i-d}(S_X)} \ar[r]^{i^X_*}  & {\o^i(X)} \ar[r]^{j^{X^*}}  &
{\o^i(X-S_X)} \ar[r] & {0} 
}
\end{align}

$j^{Y^*}(\alpha) = 0 \Rightarrow \alpha = i^Y_*(a) \Rightarrow
g^*(\alpha) = i^X_*(g^*(a)) \Rightarrow g^*(\alpha)|_{X-S_X} = 0
\Rightarrow g^*(\alpha) \in F^j(\o^i(X))$, as required. 

The proof for push-forward is similar. Let $\alpha \in F^j(\o^i(X))$. Let
$\alpha|_{X-S_X}$ is zero for codimension of $S_X \subset X$ is $d$;
$d \geqslant j$. Let $S_Y$ be the closure of $Im(f)$ in $Y$ i.e. $S_Y
\subset Y$ is push-forward of $S_X \subset X$ along $f$. Then
$dim(S_Y) \leqq dim(S_X)$ implies that codimension of $S_Y$ in $Y$ is
greater equal $j+r$. Consider the following commutative diagram 

\begin{align} 
\xymatrix{ {\o^{i-d}(S_X)} \ar[r]^{i^X_*} \ar[d]^{f_*} & {\o^i(X)}
  \ar[r]^{j^{X^*}} \ar[d]^{f_*} & {\o^i(X-S_X)} \ar[r] & {0} \\ 
{\o^{i-d+r}(S_Y)} \ar[r]^{i^Y_*}  & {\o^{i+r}(Y)} \ar[r]^{j^{Y^*}}  &
{\o^{i+r}(Y-S_Y)} \ar[r] & {0} 
}
\end{align}

Chasing the diagram in a similar way as in the pullback case
gives the result. 
\qed

\begin{Remark}\label{remark}

For $X \in Sm/k$,
$\Om(X)$ is generated as an $\Om$-module by the unit element $1 \in
\o^0(X)$ and by elements of degree greater than zero in $\Om(X)$ (\cite{levmor07},1.2.14). This
implies that if $\alpha \in 
F^j(\Om(X))$, then there exist $z_i: Z_i \ra X$ and $a_i \in \Om$ such
that $[z_i] \in \o^{\geq j}(X)$ and $\alpha = \Sigma_i a_i.[z_i]$. 

In particular, this implies that the coniveau filtration is multiplicative
(\cite{vis07}, 5.2). 

\begin{align*} F^i(\Om(X)) \cdot F^j(\Om(X)) \subset F^{i+j}(\Om(X))
\end{align*} \end{Remark}

\bl Let $X \in Sm/k$. The natural map from algebraic cobordism to Chow
rings factors 
through the coniveau filtration, i.e. \bd \o^i(X) \ra
\o^i(X)/F^{i+1}(\o^i(X)) \ra CH^i(X) \ed   \el 

\pf Given X, for any $i>0$, the subgroup $F^i(\o^i(X))$ is equal
to $\o^i(X)$. The reason is as follows. There are elements in
cobordism of the form $[f: Y \ra X, L]$ where $L$ is a line bundle over
$Y$. But the relations $\mathcal{R}_*^{sect}(X)$ implies that $[f: Y \ra
X, L]$ is equivalent to $f|_Z: Z \ra X$ for suitable choice of
codimension 1 subvariety $Z$ of $Y$. So, any element in $\o^i(X)$ is
of the form $[f: Y \ra X]$, where $dim(X) - dim(Y) = i$.  This implies
that the codimension of $Im(f)$ is at least $i$, i.e. $[f: Y \ra X] \in
F^i(\Om(X))$. 

Now, the generators of $\o^{<0} \subset \Om$ are precisely known
as maps $[V \ra {pt}]$ where $V$ is one of $\C P^n$ or Milnor
hypersurface. The pullbacks of these in $\Om(X)$ are elements of the form
$[V \times X \ra X]$. Hence the ideal $\o^{<0}.\Om(X)$ is generated
additively by the elements of the form $[V \times Y \ra Y \ra X]$, where
the first map is the projection and $[Y \ra X] \in \Om(X)$. Hence
$\o^{<0}.\Om(X) \cap \o^i(X)$ is contained in
$F^{i+1}(\Om(X))$. i.e. the kernel of the map $\o^i(X) \ra
F^i(\Om(X))/F^{i+1}(\Om(X))$ contains $\o^{<0}.\Om(X) \cap
\o^i(X)$. By remark \ref{remark}, any
element of $F^{i+1}(\Omega^i(X))$ is indeed contained in $\o^{<0}.\Om(X) \cap
\o^i(X)$. Hence $F^{i+1}(\Omega^i(X))$ is equal to $\o^{<0}.\Om(X) \cap
\o^i(X)$.  

 Hence we
have showed that the natural map from algebraic cobordism to Chow
rings factors through the coniveau filtration, i.e. \bd \o^i(X) \ra
\o^i(X)/F^{i+1}(\o^i(X)) \ra CH^i(X) \ed \qed 

We will now define algebraic cobordism of classifying spaces.

\bl\label{cod}Let $X \in Sm/k$. Then $\dfrac{\o^iX}{F^j(\o^iX)} \cong
\dfrac{\o^i(X-S_0)}{F^j(\o^i(X-S_0))}$ for a closed $S_0 \subset X$
and codimension of $S_0$ is at least j.\el 

\pf Consider the inclusion map $i: X-S_0 \hookrightarrow
X$. The localization sequence implies that the map $i^*: \o^i(X) \ra
\o^i(X-S_0) $ 
is surjective. Let $\alpha \in \o^i(X)$ and $i^*(\alpha) \in
\o^i(X-S_0)$. Let $\overline{\alpha}$ and $\overline{\i^*(\alpha)}$ be
their images in the respective quotients. 

 Suppose $\overline{\alpha} = 0$.

$\Rightarrow \alpha $ restricted to $\o^i(X-S)$ is zero for some
closed $S \subset X$ of codimension at least $j$. 

$\Rightarrow \alpha$ is zero in $\o^i(X-S_0-S)$
i.e. $\overline{\i^*(\alpha)}$ is zero. 

Suppose that $\overline{\i^*(\alpha)}$ is zero.  

$\Rightarrow$ There is  some closed $S \subset X$ of codimension at
least j such that $i^*(\alpha)$ is zero in $\o^i(X-S_0-S)$. 

$\Rightarrow \alpha$ is zero in $\o^i(X-S_0-S)$,
i.e. $\overline{\alpha}$ is zero. 

Hence $i^*:\dfrac{\o^iX}{F^j(\o^iX)} \ra
\dfrac{\o^i(X-S_0)}{F^j(\o^i(X-S_0))}$ is an isomorphism. 

\qed         

\bl\label{hom} Let $X \in Sm/k$  and $p: E \ra X$ be a
vector bundle over $X$ then, \bd \dfrac{\o^iE}{F^j(\o^iE)} \cong
\dfrac{\o^iX}{F^j(\o^iX)}\ed \el 

\pf We have that the map $p^*: \o^*(X) \ra \o^*(E)$ is an isomorphism
from the homotopy invariance. First we see that $p^*: F^j(\o^iX)
\subset F^j(\o^iE)$. For, let $x \in F^j(\o^iX)$. Let $x$ be
represented by a map $[f: Y \ra X]$ and $x|_{X-S_X}  = 0 $ for closed
$S_X \subset X$ of codimension at least $j$. Then $e = p^*(x) = [f^E:
f^*(E) \ra E]$. Let $S_E$ denote the space $E|_S$. Then $S_E$ is a
closed subset of $E$ of codimension at least $j$ in $E$.  Then
$x|_{X-S_X}  = 0 $ implies 
$e|_{E-S_E} = 0$ (from lemma \ref{pullback}). Hence $p^*: F^j(\o^iX)
\subset F^j(\o^iE)$. This implies $p^*: \dfrac{\o^iX}{F^j(\o^iX)} \ra
\dfrac{\o^iE}{F^j(\o^iE)}$ is a surjection.  

Now we will show that $p^*$ is an injection. Note that if $s_0 : X \ra E$ is
the zero section, then $s_0^*$ is the inverse map of $p^*$. This is
because as in the above set up, $f^E: f^*(E) \ra E$ is transversal to
$s_0: X \ra E$, and pulls back to $f: Y \ra X$.  

We have $e = p^*(x)$. To show an injection, we need to show if $e \in
F^j(\o^i(E))$ then $x \in F^j(\o^i(X))$. By remark \ref{remark}, we
can write  $e = \Sigma_i a_i \cdot  e_i$, where 
$a_i \in \Om$ and $e_i \in \o^{\geq j}(E)$. As $p^*$ is isomorphism,
let $e_i = p^*(x_i)$. Then $x = \Sigma a_i \cdot x_i$. As each  $x_i
\in \o^{\geq j}(X)$,  $x \in F^j(\o^i(X))$.

 \qed

\begin{Definition}\label{defn} Let $G$ be a linear algebraic group
  over a field $k$ 
  admitting resolution of singularities. Let $V_j$ be any
  representation of $G$ over $k$ such that $G$ acts freely outside a
  $G$-invariant closed suset $S_j \subset V_j$ of codimension at least
  $j$. Suppose that the geometric quotient $\frac{V_j-S_j}{G}$ exists as 
  an element of $Sm/k$.  Then we define \bd \o^i(BG) = \varprojlim_j
  \dfrac{\o^i(\frac{V_j-S_j}{G})}{F^j(\o^i(\frac{V_j-S_j}{G}))} \ed
\end{Definition}  

For a given $j$, we call a pair $(V_j,S_j)$ in the above definition as
a $j$-admissible pair over $X$.
The following theorem will show that $\o^i(BG)$ is well-defined. We
follow the technique used in \cite{tot99}.

\bt\label{independent} For $G$, $V$ and $S$ as above $ Q^j(\o^i(BG)) :=
\dfrac{\o^i(\frac{V_j-S_j}{G})}{F^j(\o^i(\frac{V_j-S_j}{G}))} $ is
independent of the choice of the representation $V_j$ and a closed
subset $S_j$ of codimension at least j \et 

\pf

Having fixed a representation, independence of the choice of $S$ is
established in lemma \ref{cod}. To obtain independence over
a choice of
representation, consider any 2 representations $V$ and $W$. Assume that $G$
acts freely outside $S_V$  of codimension at least $j$ in $V$ and
outside $S_W$ of codimension at least $j$ in $W$. Also assume
quotients $(V-S_V)/G$ and 
$(W-S_W)/G$ exists as elements of $Sm/k$. Then consider the direct sum $V \oplus
W$.  The quotient $((V-S_V) \times W)/G$ exists as an element of
$Sm/k$ being a vector 
bundle over $(V-S_V)/G$ and so also $((W-S_W) \times
V)/G$. Independence over the choice of $S$ for the representation $V \oplus W$
implies that \bd Q^j(((V-S_V) \times W)/G) \cong Q^j(((W-S_W)\times
V)/G). \ed And lemma \ref{hom} above shows that $Q^j$ of a vector
bundle is same as that of the base scheme. Hence both $V$ and $W$ have
isomorphic quotients $Q^*$ in indices less than equal to $j$,
proving independence of the choice of representation.   

\qed

\bl\label{sect} Let $X \in Sm/k$. Let $p: L \ra X$ be a line bundle
over $X$. Then 
$p^*([f: Y \ra X, f^*(L)]) = [s_0 \circ f:Y \ra L]$ where $s_0: X \ra
L$ is the zero section.\el 

\noindent{\sc Proof:} Consider the following commutative diagram
\begin{align*} 
\xymatrix{{f^*(L)} \ar[r] \ar[d]_q & L \ar[d]^p \\ Y \ar[r]_f & X. }  
\end{align*}

 We get that \bd p^*([f:Y \ra X, f^*(L)]) = [f^*: f^*(L) \ra L,
 q^*(f^*(L))].\ed 

 Let $y \in Y$ and let $q^{-1}(y)$ denote the fiber over $y$ in $f^*(L)$. Then
 the fiber over $y$ in $q^*(f^*(L))$ is precisely $q^{-1}(y) \times
 q^{-1}(y)$. Hence the diagonal map $\Delta: q^{-1}(y) \ra q^{-1}(y)
 \times q^{-1}(y) $ induces the diagonal section $s_{\Delta}: f^*(L) \ra
 q^*(f^*(L))$, transversal to the zero section. The zeros of $s_{\Delta}$ is
 precisely the closed subscheme given by the image of the zero section
 $\tilde{s_0}: Y \ra f^*(L)$. This implies \bd p^*([f: Y \ra X,
 f^*(L)]) = [s_0 \circ f:Y \ra  L]. \ed \qed 

\bl\label{chern} For $X \in Sm/k$ and a vector bundle $p: E \ra X$ of rank $n$,
$p^*(c_n(E)) = [s_0: X \ra E]$ where $s_0$ is the zero section. \el

\pf Proof is by induction on $n$. Lemma \ref{sect} implies the case $n
= 1$. Let $P(E)$ be the projectivization of $E$ and let $t: P(E) \ra X$
be the associated projective bundle map. Then we have an exact
sequence of bundles over $P(E)$
\begin{align*}
\xymatrix{
  0 \ar[r] & E_1 \ar[r]^{h_1} & t^*(E) \ar[r]^{h_2} & E_{n-1} \ar[r] & 0}
\end{align*}
where $E_1$ is the $\mathcal{O}(-1)$ bundle over $P(E)$ and rank of
$E_{n-1}$ is $n-1$. Let $q_{n-1}: E_{n-1} \ra P(E)$ and $q_{1}: E_1
\ra P(E)$ denote the bundle maps. Note that if $q: t^*(E) \ra P(E)$ is
the bundle map associated with $t^*(E)$, then $q = q_{n-1} \circ h_2$.
We have a cartesian diagram
\begin{align}\label{dig1}
  \xymatrix{ E_1 \ar[r]^{h_1} \ar[d]_{q_1} & t^*(E) \ar[d]^{h_1} \\
P(E) \ar[r]^{s_0^{n-1}} & E_{n-1}
}
\end{align}
 where $s_0^{n-1}$ is the zero section of the bundle $E_{n-1}$. Then
 \begin{align*}
   q^*(c_n(t^*(E))) =& q^*(c_1(E_1)c_{n-1}(E_{n-1})) \\
 =& h_2^*(q_{n-1}^*(c_1(E_1)c_{n-1}(E_{n-1}))) \\
=& h_2^*(s_{0_*}^{n-1}(c_1(E_1))) \text{ by induction} \\
=& h_{1_*}(q_1^*(c_1(E_1))) \text{ by \ref{dig1}} \\
=& [s_0^q: P(E) \ra t^*(E)]
 \end{align*}
where the final map $s_0^q$ is the inclusion of the zero section of
$t^*(E)$. Finally note that the pullback of $[s_0: X \ra E]$ under the map
$t^*(E) \ra E$ is  $[s_0^q: P(E) \ra t^*(E)]$. Hence we have shown
that $p^*(c_n(E)) = [s_0: X \ra E]$.
\qed

\bl For $X \in Sm/k$ and a vector bundle $p: E \ra X$, $c_i(E)
\in F^{i}(\o^{i}(X))$. \el 
\noindent{\sc Proof:} We prove this by an induction on the rank of
$E$. Let $E$ be a rank $n$ vector bundle over 
$X$. Let $s_0: X \ra E$ be the zero section. Then $p^*(c_n(E)) = [s_0:
X \ra E]$. Hence $p^*(c_n(E)) \in F^{n}(\o^{n}(E))$. By lemma \ref{hom},  $c_n(E)
\in F^{n}(\o^{n}(X))$.

For smaller $c_i$'s, we use an inductive definition of the
smaller $c_i$'s. Let $E_0$ be $E - s_0(X)$. There is a canonical rank
$n-1$ bundle $E_{n-1}'$ over $E_0$. (where the fiber over a point $(x,v)
\in E_0$ is the space orthogonal to $v$ in $E_x$). We have a 
bundle map $p: E_0 
\ra X$ and the localization sequence implies that $p^*: \o^*(X) \ra
\o^*(E_0)$ is an isomorphism in degrees less than or equal to
$n-1$. And then $c_i(E) = (p^*)^{-1}(c_i(E_{n-1}'))$ for $i < n$. Lemma
\ref{cod} implies that for $j \leqq n $,  
$\o^*(E_0)/F^j(\o^*(E_0))$ equals $\o^*(E)/F^j(\o^*(E))$ and lemma
\ref{hom} implies that $p^*$ induces an isomorphism between
$\o^*(X)/F^j(\o^*(X))$ and $\o^*(E)/F^j(\o^*(E))$. Hence by induction   
$p: E \ra X$, $c_i(E)
\in F^{i}(\o^{i}(X))$.

\qed

The Chern classes in algebraic cobordism are defined using the projective
bundle formula, but that definition is same as the one used
above. Consider $t: P(E) \ra X$, where $P(E)$ is the projectivisation of
$E$. There is an exact sequence $0 \ra O(-1)(P(E)) \ra t^*(E) \ra
E_{n-1} \ra 0$. For  the line bundle $E_1 = O(-1)(P(E))$, the space $E_0$ used
in the above proof is the total space of $L$ minus the zero section. If we
look at $t^*(E)$ as a bundle over $E_1$, i.e. pull-back of $E_{n-1}$ over
$L$, then its restriction to $E_0$ is precisely $E_{n-1}'$ used in the 
above proof. $c(t^*(E)) = c(E_1)c(E_{n-1})$ implies $c_i(E) =
(p^*)^{-1}(c_i(E_{n-1}'))$ for $i < n$ as the multiples of $c_1(E_1)$ are
mapped to zero in $\Om(E_0)$.

\section{Classical Lie Groups}\label{lie}

In this section, we assume the base field to be the complex numbers. In this
section, we will first prove a lemma which helps in finding
generators of $\Om(X)$ from that of $CH^*(X)$.  We will use this to
show that the algebraic cobordism of $BGL(n), BSL(n)$, $BSp(n)$, $BO(n)$
and $BSO(2n+1)$ respectively maps isomorphically to the complex cobordism
of these spaces.  

 In (\cite{yag05}, section 9), Yagita has considered the groups that
 satisfy $\Om_{BP}(BG) \cong BP^*(BG)$.  Here $BP^*(-)$ is
 Brown-Peterson cohomology for a prime $p$ 
 and $\Om_{BP}(-)$ is algebraic Brown-Peterson theory defined in
 section 8 of \cite{yag05}. There is a relation $\Om_{BP}(X) \cong BP^*
 \otimes_{\Om_{(p)}} \Om(X)_{(p)}$ where $\Om(X)_{(p)} := \Om(X)
 \otimes_{\Z} \Z_{(p)}$ and $\Z_{(p)}$ is the localisation of the ring
 $\Z$ at the prime $p$. Hence our calculations also implies results in
 (\cite{yag05}, section 9).  

Hu and Kriz have calculated $MGL^{*,*}(BGL(n))$ in (\cite{hukri01}, section 2).

 Firstly, we will look at the situation in $MU^*(-)$ which tells us
 what answers to expect and why. \newline

$MU^*(BGL(1))$ = $MU^*(\C P^{\infty})$ = $MU^*[[x]]$; where $x \in
MU^2(\C P^{\infty})$ is the Euler class of the canonical bundle $\gamma_1$
over $\C P^{\infty}$. For X denoting the product of $n$ copies of $\C
P^{\infty}$, \bd MU^*(X) \cong MU^*[[x_1, ... , x_n]] \ed where each $x_i \in
MU^{2}(X)$ is the pullback along the projection on the $i$th
co-ordinate of the Euler class of $\gamma_1$. Then \begin{align}
  MU^*(BGL(n)) \cong MU^*[[c_1,...,c_n]] \end{align} where each $c_i$
corresponds to the $i$th Chern class of the canonical bundle $\gamma_n$
over $BGL(n)$ which is homotopy equivalent to $Gr(n,\infty)$. For each
$i$, the pullback of $c_i$ along the embedding of torus $GL(1)^n \ra
GL(n)$ is the 
$i$th symmetric polynomial in $x_j$'s. $BSL(n)$ is homotopy equivalent
to a $GL(1)$-bundle over $BGL(n)$. This $GL(1)$-bundle is the complement of
the zero section in the total space of $det(\gamma_n)$ over $BGL(n)$. This
helps in concluding using the topological version of lemma
\ref{enot} that \begin{align}MU^*(BSL(n)) \cong
  MU^*[[c_2,...,c_n]]. \end{align}

Similarly, we have $MU^*(BSp(1))$ = $MU^*(BSL(2))$ = $MU^*[[c_2]]$
where $c_2$ is contained in $MU^{4}(BSp(1))$. And \begin{align} MU^*(BSp(n)) \cong
  MU^*[[c_2, ... , c_{2n}]] \end{align} where $c_{2i}$ is the
$2i$th Chern class of the standard representation of $Sp(n)$. If $Y$
is the product of $n$ copies of $BSL(2)$, then \bd MU^*(Y) \cong
MU^*[[y_1,...,y_n]] \ed where each $y_j \in MU^4(Y)$. The pullback of $c_{2i}$
along the diagonal inclusion $SU(2)^n \ra Sp(n)$ is the $i$th symmetric
polynomial in $y_j$'s. \newline

The $BO(n)$ case is done by Wilson in \cite{wil84}. He has shown that
$MU^*(BO(n))$ is isomorphic to the quotient of $MU^*(BGL(n))$ by the
relations $c_i = c_i^*$  where $c_i^*$ is the $i${th} Chern class of
the dual of the standard representation of $GL(n)$. These relations
are forced because the standard representation of $O(n)$ is self
dual. \begin{align} MU^*(BO(n)) \cong MU^*[[c_1,...,c_n]]/(c_1-c_1^*,  ..., c_n -
  c_n^*).   \end{align}  The expression for $c_i^*$ can be obtained as 
follows. Let $z_j = [-1](x_j)$ be the inverse of $x_j$ under the formal
group law. Then pullback of $c_i^*$ to $MU^*(X)$ is the $i${th}
symmetric polynomial 
in $z_j$'s i.e. in $[-1](x_j)$'s which can be written in terms of the
symmetric polynomials in $x_j$'s. The $i$th symmetric polynomial in
$x_j$'s is the pullback of $c_i$. This gives the expression for $c_i^*$ as a
formal power series in $c_i$'s with coefficients in $MU^*$. The
topological version of lemma \ref{enot} implies that
\begin{align}
  MU*(BO(2n+1)) \cong MU*[[c_2,\cdots,c_n]]/(c_2 - c_2^*,\cdots,c_n-c_n^*).
\end{align}

Now we will start our calculations.

\bl \label{gen} Let $X \in Sm/k$. If $CH^*X$ is generated as an
abelian group by some 
elements $x_1,x_2,\cdots$ and if $y_1,y_2,\cdots \in \Om X$ map to
$x_i \in CH^*X$ then $\Om X$ is generated as a $\Om$-module by
$y_1,y_2,\cdots$.\el 

\pf    We have a natural map $\Om X \ra CH^*X$
and there is an isomorphism  $CH^*X = \Om X
\otimes_{\Om}\Z$. This implies that the kernel of the map $\Om X \ra CH^*X$
is $\o^{<0} \cdot \Om(X)$. We prove the lemma by induction on decreasing degree of
$y \in \Om X$. Let the image of $y$ in $CH^*X$ be $x = \Sigma
n_ix_i$. Then $y-\Sigma n_iy_i$ is mapped to zero in the Chow ring. Hence
$y - \Sigma n_iy_i = \Sigma u_jz_j$ where $u_j \in
\o^{<0}$. Hence each of the $z_j$ have the degree higher than
the degree of $y$. As
 cobordism is zero in the degrees higher than dimension of $X$, we
are through by induction. 
\qed \newline

 We use this lemma to calculate the algebraic cobordism of $BGL(n), BSL(n)$,
 $Sp(n)$, $BO(n)$ and $BSO(2n+1)$. The Chow rings of these spaces have
 been calculated in \cite{tot99} and then also in \cite{molvis06} using
 an unified approach termed as the `stratification'. 

$CH^*(BGL(n)) \cong \Z[c_1,c_2,...,c_n]$; where each $c_i$ is the $i$th
Chern class of the canonical $n$-bundle $\gamma_n$ over $BGL(n)$ in the Chow
ring. Since the canonical map $\o^*X \ra CH^*X$ is functorial, for any
vector bundle $E$ over $X$, it maps the Chern classes of $E$ in $\o^*(X)$
to the Chern classes of $E$ in $CH^*X$. Note that Levine-Morel defined the
Chern classes as maps $c_i: \o^*X \ra \o^{*+i}X$ (\cite{levmor07},
4.1.15). Here by the Chern class as element in the algebraic cobordism, we
take $c_i(id)$, and continue to denote it by $c_i$ where $id \in \Om(-)$ is
the identity element. 

Hence we know that $\o^*(BGL(n))$ is generated as an $\o^*$-module by the
monomials in $c_1,c_2,...,c_n$; the Chern classes of $\gamma_n$ over
$Gr(n, \infty)$. Now consider the natural map \bd \o^*(BGL(n)) \ra
MU^*(BGL(n)) \cong MU^*[[c_1,c_2,...,c_n]].\ed We know that the natural
map $\o^*(pt) \ra MU^*(pt)$ is an isomorphism. Hence any non-zero
polynomial or power series relation between $c_i$'s with coefficients
in $\o^*$ will be mapped to non-zero relation between $c_i$'s with
coefficients in $MU^*$. But since there is no such relation in the complex
cobordism ring of $BGL(n)$, there is no such relation in algebraic
cobordism too. Hence 
\begin{align} \o^*(BGL(n)) \cong \o^*[[c_1,c_2,...,c_n]]. \end{align} 

The exact similar arguments, i.e. obtaining the Chern classes as the generators
using the Chow ring and using the complex cobordism to show there is no
relation between Chern classes implies that \begin{align} \o^*(BSL(n)) \cong&
\o^*[[c_2,...,c_n]]  \\ \o^*(BSp(2n)) \cong& \o^*[[c_2, c_4,...,
c_{2n}]]. \end{align}     

We can use this method for $BO(n)$ and $BSO(2n+1)$ too using the
calculations of $MU^*BO(n)$ done in [13].  

We have $CH^*BO(n) \cong \Z[c_1,...,c_n]/(2c_i=0 ; i$ odd). This implies that
$\o^*(BO(n))$ as a module over $\Om$ is generated by the
monomials$c_1,...,c_n$; Chern classes in 
algebraic cobordism of standard representation of $BO(n)$. Now we
obtain the relations among these Chern classes. Let $c_i^*$ be the
Chern classes for the dual of standard representation of $GL(n)$.  As the
 standard representation $O(n) \ra GL(n)$ is self-dual, pullbacks of
 $c_i$ and $c_i^*$ are equal. So, we get $c_i = c_i^*$ in algebraic
 cobordism too. 

Now consider the natural map $\o^*(BO(n)) \ra MU^*(BO(n))$ and the
isomorphism $\o^* \ra MU^*$. Absence of any polynomial or power series
relation besides $c_i=c_i^*$ in $MU^*(BO(n))$ implies that there are no
more relations between $c_i$'s other than $c_i = c_i^*$ in algebraic
cobordism too. Hence \begin{align} \o^*(BO(n)) \cong
\dfrac{\o^*[[c_1,...,c_n]]}{<c_i= c_i^*>} \end{align} 

Similarly we get that \begin{align} \o^*(BSO(2n+1)) \cong
\dfrac{\o^*[[c_2,c_3,...,c_{2n+1}]]}{<c_i= c_i^*>}. \end{align}

\section{Products of finite cyclic groups}

In this section we continue to assume that the base field is the
complex numbers. We will calculate $\Om(BG)$ when $G$ is a product of finite
cyclic groups. $CH^*(BG)$ for such $G$ is calculated in \cite{tot99} and
$MU^*(BG)$ is calculated in \cite{lan70}.  As
said before, Yagita has calculated $\Om_{BP}(BG)$ for $G$ a product of
finite abelian $p$-groups and our calculations imply results in
(\cite{yag05}, section 9).    

In this section,  we will follow the approach used in \cite{lan70}
and in  the process show that a version of Kunneth formula for Chow rings
proved by Totaro in \cite{tot99} also holds in algebraic cobordism.

\bl\label{enot} Let $X \in Sm/k$ and $p: E \ra X$ be a
vector bundle of rank r. Let 
$E_0$ be a complement of the zero section of E. Then the pullback homomorphism
$\Om(X) \ra \Om(E_0)$ is always surjective and its kernel is generated
by the top Chern class $c_r(E)$. \el

\pf This follows from the localization sequence, Let $s_0: X \ra E$ be the
zero section and let $j: E_0 \ra E$ be the inclusion map. Let $p_0:
E_0 \ra X$ be the restriction of $p$. The localisation sequence
corresponding to $j$ is

\begin{align*}
\xymatrix {{\Om(X)} \ar[r]^{s_{0_*}} & {\o^{*+r}(E)} \ar[r]^{j^*} &
  {\o^{*+r}(E_0)} \ar[r] & {0.}}  
\end{align*}

By the homotopy invariance, we have an isomorphism $p^*:\Om(X) \ra
\o^{*}(E)$. Hence by the localisation sequence, the  pullback along $p_0$,
denoted by $p_0^*: \Om(X) \ra \Om(E_0)$ is always surjective. To
show that the
kernel of this surjection is generated by $c_r(E)$, we need to show
$s_0^*(x) = p^*(c_r(E)\cdot x)$ for $x \in \Om(X)$. Lemmas \ref{sect} and
\ref{chern} show exactly this. Hence proved.

\qed

$\Omega^*(\mathbb{C}P^\infty)$ = $\Omega^*[[x]]$, where $x$ is the first
Chern class of $\mathcal{O}(1)$ bundle over $\mathbb{C}P^\infty$,
i.e. $x = c_1(\mathcal{O}(1))(id)$. By $[n](x)$, we mean the addition of
$x$ with itself $n$ times, under the formal group law  in
$\Omega^*(\mathbb{C}P^\infty)$, $F(a,b) = a + b +
\sum_{i,j>0} a_{ij}a^ib^j$. Then
\begin{align*}[2](x) = F(x,x) = 2x + \sum_{i,j>0} a_{ij}x^{i+j} =
  2x + a_{11}x^2 + 2a_{12}x^3 + \cdots. \end{align*} 

\begin{center}$[3](x) = F(x, [2](x)) = 3x + 3a_{11}x^2 + (a_{11}^2 +
  8a_{12}) x^3 + \cdots$\end{center} and $[n](x) = F(x,[n-1](x))
= nx +$ higher order terms.   
  
\begin{Lemma} $\Omega^*(B\mathbb{Z}_n) = \Omega^*[[x]]/([n](x))$. \end{Lemma}

\noindent{\sc Proof:} If $E$ is the total space of $\mathcal{O}(n)$
over $\mathbb{C}P^\infty$ and $s_0: \mathbb{C}P^\infty \rightarrow E$
is the zero section then,  $B\mathbb{Z}_n$ is $E - Im(s_0)$.

The localization sequence implies that \begin{displaymath}
\om(\mathbb{C}P^\infty) \rightarrow^{f} \om(E) \rightarrow
\om(B\mathbb{Z}_n) \rightarrow 0 \end{displaymath}

Considering the isomorphism $\om(E) \cong
\om(\mathbb{C}P^\infty)$, we get that, the map $f$ is
 the multiplication by $c_1(\mathcal{O}(n))$.

From the formal group law, $c_1(\mathcal{O}(n)) = [n](x)$. Hence
\begin{displaymath} \Omega^{*}(B\mathbb{Z}_n) =
\Omega^*[[x]]/([n](x)).\end{displaymath}

\qed

\bl\label{kunneth} Let $X$ and $Y$ be in $Sch/k$, and that
$Y$ can be partitioned into open subsets of an affine space, then the
cup product map $\Omega^*(X) \otimes_\Om \Om(Y) \ra \Om(X \times Y)$
is surjective. 
\el
\noindent{\sc Proof:} Using the localization sequence for an open
subset $U$ of $\A^n$, if $Z$ is the complement of $U$,  $\Om(Z) \ra
\Om(A) \ra \Om(U) \ra 0$ implies that $\Om(A) \ra \Om(U)$ is a
surjection. The homotopy invariance property says that $\Om(X)
\otimes_{\Om} \Om(\A^n) \cong \Om(X \times \A^n)$. As $X \times U$ is
an open subset of $X \times \A^n$, the localization sequence implies
that $\Om(X \times \A^n) \ra\Om( X \times U)$ is also
surjection. Hence the commutative diagram 

\begin{align}\xymatrix{
{\Om(X) \otimes_\Om \Om(\A^n)} \ar[r]^{\cong} \ar[d] & {\Om(X \times
  \A^n)} \ar[d] \\ 
{\Om(X) \otimes_\Om \Om(U)} \ar[r] & {\Om(X \times U)}
}\end{align}
 implies that ${\Om(X) \otimes_\Om \Om(U)} \ra  {\Om(X \times U)}$ is
 surjective.  
 
 Suppose that $Y$ can be partitioned into open subsets of affine
 space. Let $U$ be an open subset among that partition and $Z$ be its
 complement in $Y$.  

\begin{align}\xymatrix{
{\Om(X) \otimes_\Om \Om(Z)}  \ar[r] \ar[d]  & {\Om(X) \otimes_\Om
  \Om(Y)} \ar[r] \ar[d] & {\Om(X) \otimes_\Om \Om(U)} \ar[d] \ar[r] &
0\\ 
{\Om(X \times Z)} \ar[r] & {\Om(X \times Y)} \ar[r] & {\Om(X \times
  U)} \ar[r] & 0 
}\end{align}
 
The above diagram implies that $\Om(X)\otimes\Om(Y) \ra \Om(X \times Y)$
is surjective if both $\Om(X)\otimes\Om(Z) \ra \Om(X \times Z)$ and
$\Om(X)\otimes\Om(U) \ra \Om(X \times U)$ are surjective. Hence we
obtain the result by induction on the number of open subsets in the
partition.  

\qed

$B \Z /n$ can be taken to be the complement of zero section in $O(n)$
bundle on $\C 
P^{\infty}$. It can be successively approximated by the complement of
the zero
section in the $O(n)$ bundle on the finite dimensional 1-Grassmannians,
i.e. complex projective spaces (\cite{tot99}, 1.4). Since $\C P^n$'s
can be partitioned into affine spaces, each of the approximations of
$B \Z/n$ can be partitioned into open subsets of affine
spaces. \newline 

 If $X$ is an infinite CW complex with finite skeletons $X_n$, then
 these skeletons induce a filtration on $MU^*(X)$. i.e. we say $x \in
 MU^*(X)$ is in the $n^{th}$ level of filtration if if its zero when
 restricted to $X_n$. From \cite{lan70},  if $MU^*(X)$ has no
 elements of infinite filtrations then \begin{align} MU^*(B\Z/n
   \times X) \cong& MU^*(B\Z/n)\otimes_{MU^*}MU^*(X) ; \\ MU^*(B\Z/n
   \times B\Z/m) \cong& MU^*(B\Z/n)\otimes_{MU^*}MU^*(B\Z/m).
\end{align} 

\bl  $\Om(B\Z/m) \otimes_\Om \Om(B\Z/n)  \cong \Om(B(\Z/m \times \Z/n)) $.
\el
\noindent{\sc Proof:} For a complex smooth variety $X$, we have a
natural map $\Om(X) \ra MU^*(X(\C))$. This map is an isomorphism for
$X = B\Z/n$ and $Y = B\Z/m$.  
Also note that \bd MU^*(X) \otimes_{MU^*} MU^*(Y) \ra MU^*(X \times Y) \ed
is an isomorphism. Now consider the diagram  

\begin{align*}\xymatrix{
{\Om(X) \otimes_{\Om} \Om(Y)} \ar[r] \ar[d] & {\Om(X \times Y)} \ar[d] \\
{MU^*(X) \otimes _{MU^*} MU^*(Y)} \ar[r] & {MU^*(X \times Y)} 
}\end{align*}

The left vertical map is an injection and the lower horizontal map is
an isomorphism in case $X$ and $Y$ are classifying spaces of cyclic
groups. This implies that the upper horizontal map is an
injection. Surjection is proved by the lemma \ref{kunneth}. Hence
$\Om(B\Z/m) \otimes_\Om \Om(B\Z/n)  \cong \Om(B(\Z/m \times \Z/n)) $ 

\qed

\noindent Hence now we have,
\bt If $G = \Z/n_1 \times \Z/n_2 ... \times \Z/n_r$ then $\Om(BG) =
\dfrac{\Om[[x_1, x_2, ..., x_r]]}{([n_1](x_1),..., [n_r](x_r) ) } $
where $x_i \in \o^1(BG)$ for all $i$. 
\et

\section{Cellular spaces} 

  Let $X \in Sch/k$ and suppose there exists  a filtration $\phi = Z_0
  \subset Z_1 \subset 
 ... \subset Z_n = X$ where $U_i = Z_{i+1} - Z_i$ is an  affine
 space $\A^{N_i}_k$ (a cell). We call such  $X$ a cellular space. (\cite{levmor07},
 5.2.11) shows that $\om(X) \ra \bigoplus_{U_i} \o_*({U_i})$ induced
 by restriction maps is a surjection over any base field allowing
 resolution of singularities. We will show its an isomorphism over the
 complex numbers
 using the natural map to the complex cobordism. As expected, we get
 \bd\om(X) \cong MU_*(X).\ed The statement is known for a long time now,
 but for the sake of completeness we provide the proof. For the next
 theorem, we
 assume that the base field is the complex numbers.

\bt For above $X$, $\o_r(X)$ is isomorphic to  $\bigoplus_{U \in I}
\o_r({U})$, where $I$ is set of cells of $X$. In fact, the map
$\om(X) \ra MU_{2*}(X)$ is an isomorphism.\et 
 
\noindent{\sc Proof:} When there is just one cell, statement follows
trivially. We now proceed by induction on the number of cells. Suppose
the 
statement holds for $Z$ and we want to add a cell $U$ of dimension $n$
to get $X$, i.e. $Z \hookrightarrow X \hookleftarrow U$. Let $\bar{U}$
be closure of $U$ in $X$. And $\bar{U}'$ be a resolution of singularities
of $\bar{U}$.  Consider the diagram

\begin{align}
\xymatrix {Z \ar@{^{(}->}[r] &  X & U \ar@{_{(}->}[l] \\
{\bar{U}'} \ar[r] \ar[drr]_{\tilde{f}} & {\bar{U}} \ar[u]^{i}  & U
\ar[l]^j \ar[u]_{\cong} \ar[d]^f \\ 
& & (Pt)
}
\end{align}

 This gives the following commutative diagram on cobordism groups,

\begin{align}\xymatrix {{\o_r(Z)} \ar[r] & {\o_r(X)} \ar[r] & {\o_r(U)} \ar[r] & 0 \\
{\o_r(\bar{U}')} \ar[r] & {\o_r(\bar{U})} \ar[u]^{i_*} \ar[r]^{j^*} &
{\o_r(U)} \ar[u]_{\cong}  \\ 
& & \o_{r-n}(Pt) \ar[ull]^{\tilde{f}^*} \ar[u]_{\cong}^{f^*}
}
\end{align}
 
Hence $i_* \circ \tilde{f}^*$ gives a section of $\o_r(U)$ in
$\o_r(X)$. Which implies $\o_r(X) \cong \o_r({U}) \oplus
Im(\o_r(Z))$. To prove that $\o_r(Z)$ injects into $\o_r(X)$ we use
the natural isomorphism from $\om$ to $MU_*$. 

As the natural map $\om(\A^n) \ra MU_*(\A^n)$ is an isomorphism, by
induction, we get that the map $\om(Z) \ra MU_*(Z)$ is an
isomorphism. We have a commutative diagram, 
 
\begin{align}
\xymatrix{  & {\o_r(Z)} \ar[r] \ar[d] & {\o_r(X)} \ar[r] \ar[d] &
  {\o_r(U)} \ar[r] \ar[d] & 0 \\ 
MU_{2r+1}(U) \ar[r] & MU_{2r}(Z) \ar[r] & MU_{2r}(X) \ar[r] &
MU_{2r}(U) \ar[r] & 0 
}
\end{align}

Since $U$ is a complex affine space $\A^{N}_{\C}$, $MU_{2i+1}(U)$ is zero for all
$i$. Hence we have that the map $MU_{2r}(Z) \ra MU_{2r}(X)$ is
injective. Hence we 
get that the map $\o_r(Z) \ra \o_r(X)$ is also injective. Hence \bd \om(X) \cong
\bigoplus_{U \in I} \om({U}).\ed In particular $\om(X) \ra MU_*(X)$ is
an isomorphism.   
\qed 

Note that the isomorphism in the above theorem is not unique and depends
on the choice of resolution of singularities $\bar{U}'$.  

We can apply this to projective homogeneous varieties which have
the Schubert cell decomposition. In particular we can apply this to
Grassmanians.

\section{Equivariant cobordism} 

 Let $G$ be a linear algebraic group. By a $G$-scheme, we mean  an
 element of $Sch/k$ equipped with an action of $G$. In this section we
 define $G$-equivariant algebraic cobordism similar 
 to $G$-equivariant 
Chow ring defined by Edidin and Graham in \cite{edigra98}. We show
that $G$-equivariant algebraic cobordism  satisfies the localization
exact sequence.   

Let $G$ by any linear algebraic group $G$ and a let $X$ be a smooth
$G$-scheme. We say that an action of $G$ on $X$, denoted by say
$\sigma$, is linearized if there exist a line bundle
\begin{align*}
  \pi: L \ra X
\end{align*}

over $X$ such that, the action $\sigma$ of $G$ on $X$ lifts to an action on
$L$ (\cite{mum94}, Definition
1.6). In partiucular, we have a commutative diagram
\begin{align*}
  \xymatrix{
G \times L \ar[d]_{\mathbb{1}_G \times \pi} \ar[r]^{\Sigma} & L
\ar[d]^{\pi}\\
G \times X \ar[r]^{\sigma} & X.
}
\end{align*}
 
Let $(V_j,S_j)$ be a $j$-admissible pair as in definition \ref{defn}. Let $U$ be $V_j-S_j$. Whenever the action of $G$ is linearized for a line bundle which is relatively
ample for the projection $X \times U \ra U$ and $U/G$ is known to be
quasi-projective, then the quotient $(X \times U)/G$ is
quasi-projective (\cite{mum94}, proposition 7.1). From now on, we will
take this as a definition for $G$-linearized action. Furthermore, if
both $X$ and $U$ are smooth, then the quotient $(X \times U)/G$ is
also smooth.

\begin{Definition}\label{equivariant} For a linear algebraic group $G$ and a smooth
  scheme $X$ equipped with a $G$-linearized action
 and
any $j$-admissible pair $(V_j,S_j)$ as in definition \ref{defn}, we define
\begin{align} \o^i_G(X) = 
\varprojlim_j \dfrac{\o^i(\frac{X \times
    (V_j-S_j)}{G})}{F^j(\o^i(\frac{X \times
    (V_j-S_j)}{G}))}. \end{align}   \end{Definition}

\bl Let $G$ be a linear algebraic group with normal subgroup $H$ and let $X$
be a smooth $G$-scheme with a $G$-linearized action. Suppose the
action of H on X is free and the quotient 
scheme $X/H$ exists. Then \bd \o^*_G(X) \cong \o^*_{G/H}(X/H).\ed\el 

\pf Let $(V_j,S_j)$ be a $j$-admissible pair over $G$ for a given $j$.

\begin{align*}
\o^i_G(X) =& \varprojlim_j \dfrac{\o^i(\frac{X \times
    (V_j-S_j)}{G})}{F^j(\o^i(\frac{X \times (V_j-S_j)}{G}))} \\
=&  \varprojlim_j \dfrac{\o^i(\frac{(X \times
    (V_j-S_j))/H}{G/H})}{F^j(\o^i(\frac{(X \times
    (V_j-S_j))/H}{G/H}))} \\
=&  \varprojlim_j \dfrac{\o^i(\frac{(X/H) \times
    (V_j - S_j)}{G/H})}{F^j(\o^i(\frac{(X/H) \times (V_j-S_j)}{G/H}))}
\end{align*} 
where the last equality follows because $(X \times
    (V_j-S_j))/H$ is isomorphic to $(X/H) \times (V_j-S_j)$ as the
    action of $H$ on $X$ is free.

The final expression is by definition isomorphic to $\o^i_{G/H}(X/H)$, which is required. 

\qed

Now we have to check that the $G$-equivariant localization sequence also holds.

\bt\label{equilocalization} Let $X$ be a smooth scheme with a $G$-linearized action. Let $Z$
be its closed $G$-subscheme of 
codimension $r$ and let $U$ be the complement $Z$ in $X$. Then \bd
\o^i_G(Z) \ra 
\o^{i+r}_G(X) \ra \o^{i+r}_G(U) \ra 0 . \ed \et  

\pf  Let $(V,S)$ be a $j$-admissible pair over $G$ for a given $j$. Let
$X'$ denote the 
scheme $\frac{X \times 
  (V-S)}{G}$ and similarly denote by $Z'$ and $U'$ the subschemes of
$X'$ corresponding to $Z$ and $U$. To prove the lemma, it is enough to
show that for each $j$,\begin{align}
\dfrac{\o^i(Z')}{F^j(\o^i(Z))} \ra
\dfrac{\o^{i+r}(X')}{F^j(\o^{i+r}(X'))} \ra
\dfrac{\o^{i+r}(U')}{F^j(\o^{i+r}(U'))} \ra 0. \end{align} This follows from
analysing the following commuative diagram  

\begin{align}\xymatrix { {\o^i(Z')} \ar[r]^{i_*} \ar[d] &
    {\o^{i+r}(X')} \ar[r]^{j^*} \ar[d] & {\o^{i+r}(U')} \ar[r] \ar[d]
    & {0} \\  
{\dfrac{\o^i(Z')}{F^j(\o^i(Z))}} \ar[r]^{i'_*} &
{\dfrac{\o^{i+r}(X')}{F^j(\o^{i+r}(X'))}} \ar[r]^{j'^*} &
{\dfrac{\o^{i+r}(U')}{F^j(\o^{i+r}(U'))}} \ar[r] & {0} \\ 
} \end{align}  

$Im(i_*) = Ker(j^*)$ implies $Im(i'_*) \subset Ker(j'^*)$. 

Now we will show the containment the other way. Let $x \in
\o^{i+r}(X')$. Let its image in the quotient
$\dfrac{\o^{i+r}(X')}{F^j(\o^{i+r}(X'))}$ be $\overline{x}$ and
$j^*(x) 
= u$. Suppose $j'^*$ maps $\overline{x}$ to zero, then we need to show
$\overline{x}$ is the image of some $\overline{z} \in \dfrac{\o^i(Z')}{F^j(\o^i(Z))}$, where $z \in
\o^i(Z') $.  

The idea is to find another cobordism cycle $x' \in F^j(\o^{i+r}(X'))$
s.t. $j^*(x') = u$. So that we can write $x - x' = i_*(z)$, and then
taking the images in the lower row gives $\overline{x} =
i_*'(\overline{z})$. We find $x'$ as follows. 

$j'^*$ maps $\overline{x}$ to zero implies that $u \in
F^j(\o^{i+r}(U'))$. Hence there is a subset $S$ of $U'$ s.t., it is closed
in $U'$, its codimension in $U'$ denoted by $r_S$ is greater equal
$j$ and $u|_{U-S} = 0$. This means $u = i^S_*(s)$ for some $s \in \Om(S)$. Let
$\overline{S}$ be the closure of $S$ in $X'$. We have a cartesian square
\begin{align}\xymatrix {{\overline{S}} \ar[r]^{i^{\overline{S}}} &
    {X'} \\ {S} \ar[r]_{i^S} \ar[u]_f & {U'} \ar[u]^j} \end{align}
where all the arrows are inclusions and it gives a commutative diagram
of the respective cobordism groups
\begin{align}\xymatrix {\o^{i+r-r_S}({\overline{S})}
    \ar[r]^{i^{\overline{S}}_*} \ar[d]^{f^*} &
    {\o^{i+r}(X')} \ar[d]_{j^*} \\ {\o^{i+r-r_S}(S)} \ar[r]_{i^S_*}  &
    {\o^{i+r}(U')}.} 
\end{align}

There exists $\overline{S}$ such that, $s = f^*(\overline{s})$ as
$f^*$ is a surjection by the localization 
sequence. And we set $x' = i^{\overline{S}}_*(\overline{s})$. The above
commutative diagram implies that \bd j^*(x') = i_*^S(f^*(\overline{s})) =
u.\ed
Hence we have found a correct $x'$ we needed and hence proved.

\qed

\section{The quaternion group}

 In this section, we again assume that the base field is the complex
 numbers. Let $Q = \{ \pm 1, \pm i, \pm j, \pm k\}$ be the group of 
 quaternions. In this section we will calculate $\o^*(BQ)$.

The group $Q$ has 4 one-dimensional representations $\mathds{1},
 \varrho_i, \varrho_j, \varrho_k$. $\varrho_i$ denotes the
 homomorphism $Q \ra \{\pm 1\} \subset \C$, which maps $i$ to $-1$ and
 $j,k$ to $1$. It also has a unique 2-dimensional irreducible complex
 representation $\rho$. Let $L_i, L_j$ and $L_k$ respectively denote the
 line bundles over $BQ$ corresponding to $\varrho_i, \varrho_j,
 \varrho_k$. And let $V$ be the rank-2 vector bundle corresponding to
 $\rho$. Relations between these representations are given in
 \cite{ati61}. Written in terms of corresponding vector bundles, these
 are  

\begin{align}
\label{r1}L_i \otimes L_j \cong& L_k.\\
\label{r2}L_i \otimes L_i \cong& \mathds{1} \text{ and similar for } j \text{
  and } k.\\ 
\label{r3}L_i \otimes V \cong& V \text{ and similar for } j \text{ and } k.\\
\label{r4}L_i \oplus L_j \oplus L_k \oplus \mathds{1} \cong&  V \otimes V.\\
\label{r5}det(V) \cong& \mathds{1}.
\end{align}

As shown in \cite{ati61}, $H^*(BQ, \Z) $ is generated by $c_1(L_i) = x'$,
$c_1(L_j) = y'$ and $c_2(V) = z'$. The above mentioned relations in
representations force the relations on these generators. We get for
$x',y' \in H^2(BQ, \Z)$ and $z' \in H^4(BQ, \Z)$, \begin{align} H^*(BQ, \Z) =
\Z(x',y',z')/<2x', 2y', 8z', x'^2, y'^2, x'y'-4z'>.\end{align} 

 


We will first  show \bd CH^*(BQ) \cong
H^{2*}(BQ, \Z). \ed The argument was outlined in (\cite{gui05}, 1.1.3).

\bl $ CH^*(BQ) \cong
H^{2*}(BQ, \Z).  $ \el

\pf
The natural map $CH^*(BQ) \ra H^{2*}(BQ, \Z)$ is an isomorphism for $*
= 1,2$ by (\cite{tot99}, 3.2). Consider the following commutative
diagram of the 
localisation sequences for the representation $V$.

\begin{align*}
\xymatrix{ CH^i_Q(\{0\}) \ar[r]^{\cdot c_2(V)} \ar[d] & CH^{i+2}_Q(\A^2)
  \ar[r] \ar[d]   & CH^{i+2}_Q(\A^2-\{0\}) \ar[r] \ar[d] & 0 \\
H^{2i}_Q(\{0\}) \ar[r]^{\cdot c_2^H(V)}  & H^{2i+4}_Q(\A^2) \ar[r] &
H^{2i+4}_Q(\A^2-\{0\}) \ar[r] & 0.
}  
\end{align*}

Note that $Q$ acts freely on $\A^2 - \{0\}$, and $(\A^2-\{0\})/Q$ has the
dimension 2. Hence we get that for $i \geq 1$, the map given by
multiplication by 
$c_2(V)$ is a surjective map from $CH^i(BQ) \cong CH^i_Q(\{0\}) $ to
$CH^{i+2}_Q(\{0\})$. We know that the multiplication by $c_2^H(V)$ is
always an isomorphism from $H^{2i}_Q(\{0\})$ to
$H^{2i+4}_Q(\{0\})$. Hence the commutative diagram reduces to
\begin{align*}
  \xymatrix{ CH^i_Q(\{0\}) \ar[r]^{\cdot c_2(V)} \ar[d] & CH^{i+2}_Q(\{0\})
  \ar[r] \ar[d]   &  0 \\
H^{2i}_Q(\{0\}) \ar[r]_{\cdot c_2^H(V)}^{\cong}  & H^{2i+4}_Q(\{0\}) \ar[r] & 0.
} 
\end{align*}

As we already know that the natural map $CH^*(BQ) \ra H^{2*}(BQ, \Z)$
is an isomorphism for $* = 1,2$, we have the result by induction.

\qed

Hence the lemma \ref{gen} implies that $\Om(BQ)$ is generated as an
$\Om$-module by $c_1(L_i)=x$, $c_1(L_j)=y$ and $c_2(V)=z$. What now
remains is to find all the relations on these generators.  
%


As before, let $F(a,b) = a + b + \Sigma_{i>0,j>0} a_{ij}a^ib^j$ be the formal
group law of the complex cobordism. Whenever convenient, we will
denote $F(a,b)$ by $a +_F b$. The idea is to understand $MU^*(BQ)$
first, and then show that $\Om(BQ) \ra MU^*(BQ)$ is an
isomorphism. Mesnaoui in \cite{mes90} has calculated $MU^*(BQ)$
explicitly. He has shown that there is an ideal $I$ generated by
certain six homogeneous power series in $MU^*[[c_1(L_i), c_1(L_j),
c_2(V)]]$ such that \begin{align} MU^*(BQ) = MU^*[[c_1(L_i),
  c_1(L_j), c_2(V)]]/I.   \end{align}

We will make a standard abuse of the splitting principle to explain
what these six 
power series are in a different way. The idea is following. Given a base
space $X$, and a rank $n$ vector bundle $p: E \ra X$, the projective
bundle formula implies that $MU^*(X)$ injects into $MU^*(F(E))$, for the
flag bundle of $E$, $q_1: F(E) \ra X $. If $q: P(E) \ra X$ is the 
projectivasation, then there is an exact sequence $0 \ra L_1 \ra
q_1^*(E) \ra E_{n-1} \ra 0$ for line bundle $L_1$. So $c(q_1^*(E)) =
c(L_1)c(E_{n-1})$. Repeating the argument for $E_{n-1}$ and so on, we can
claim that there are line bundles $L_i$ over $F(E)$ such that $q^*(E)$
splits as a sum $L_i$'s.  \newline

Using this for $V$ over $BQ$, there is a space $p:X \ra BQ$ such that
$MU^*(BQ)$ injects in $MU^*(X)$ and $c(p^*(V)) = c(M)c(N)$, for some line
bundles $M$ and $N$ over $X$, where $c(-)$ is the total Chern class
operator. Let $a$ and $b$ respectively denote the
pullbacks of $c_1(L_i)$ and $c_1(L_j)$ in $X$ and $m$ and $n$
respectively denote the first Chern classes of $M$ and $N$. We will now
use the formal group law to find the relations in cobordism from the relations
in representations.

We will start with relation \ref{r5}. Here, it is interesting to
note that for a vector bundle $E$, the equation $c_1(det(E)) = c_1(E)$,
which is true in $H^*(-,\Z)$ and $CH^*(-)$ does not hold in the
cobordism theory
in general. If $E$ splits as a sum of line bundles $\oplus_i L_i$, then
$det(E) = \otimes_i L_i$. Because $H^*(-, \Z)$ and $CH^*(-)$ have the
additive formal group law, in those theories we get \begin{align}
  c_1(\otimes_i L_i) = c_1(\oplus L_i). \end{align} But the formal
group law in the cobordism theory, $F(-,-)$, is not additive and $det(V) =
\mathds{1}$ implies that \begin{align} F(m,n) = 0 \text{ i.e. } m +_F
  n = 0.  \end{align}  

We can transform the equation $F(m,n) = 0$ completely algebraically by
using the relations among $a_{ij}$'s to make it look like $m+n =
P(mn)$ where $P(-)$ is a formal power series with  coefficients in
$MU^*$. This is done as follows. 
 
\bl\label{rel} Let $MU^*(BGL(1)) = MU^*[[c_1]]$. If $[-1](c_1)$
denotes the inverse 
of $c_1$ under the formal group law, then there exist a formal power
series $P(-) \in MU^*[[c_1]]$ such that $c_1+[-1](c_1) = P(c_1 \cdot [-1](c_1))$
\el 

\pf Let $f: GL(1) \ra SL(2)$ be a map given by $x \mapsto
\begin{pmatrix} x & 0 \\ 0 & x^{-1} \end{pmatrix}$. Let $e_1$ and
$e_2$ be the Chern classes of the standard representation of $SL(2)$. Then
$MU^*(BSL(2)) \cong MU^*[[e_2]]$, and let $e_1 = P(e_2)$. As the pullbacks of
$e_1$ and $e_2$ are $(c_1 + [-1](c_1))$ and $(c_1 \cdot [-1](c_1))$
respectively, we have the result.

\qed

Again note here that $e_1 \in MU^*(BSL(2))$ is not zero. But a fact
that all the coefficients of $P(-)$ are actually in $MU^{\leq 0} \cong
\o^{\leq 0}$ implies that the image of $e_1$ in $H^*(BSL(2))$ is zero
and the image of $e_1$ in
$CH^*(BSL(2))$ is zero. Hence we have \begin{align} m+n = P(mn).
\end{align} Looking at its image in $MU^*(BQ)$, it implies
\begin{align}c_1(V)  = P(c_2(V)). \end{align} 

The relation \ref{r3} implies that $c_1(L_i \otimes V) = c_1(V)$ and $c_2(L_i
\otimes V) = c_2(V)$. Which translates as
\begin{align}F(a,m) + F(a,n) = m+n \text{ and } F(a,m)F(a,n) = mn.\end{align}
\begin{align}F(b,m) + F(b,n) = m+n \text{ and } F(b,m)F(b,n) = mn. \end{align}

The existence of $P(-)$ already implies that $P(F(a,m)F(a,n)) = P(mn)$, so 
relations $F(a,m) + F(a,n) = m+n$ and $F(b,m) + F(b,n) = m+n$ are
already implied by $F(a,m)F(a,n) = mn$ and $F(b,m)F(b,n) = mn$
respectively. Hence the third relation in representations gives only 2
independent  
relations in cobordism. \begin{align}\label{1st} F(a,m)F(a,n) =& mn.  \\
  \label{2nd} F(b,m)F(b,n)  =& mn.  \end{align}      

The relation \ref{r2} implies that  \begin{align}\label{3rd}F(a,a)
  =& 0. \\ \label{4th} F(b,b)   =& 0 .\end{align}      

Relation \ref{r4} gives four relations in the cobordism by comparing
all 4 Chern classes of the two bundles. Comparing the first Chern
class gives 
\begin{align}\label{5th1}a + b + c = [2](m) + [2](n).\end{align} Here
by $c$, we mean pullback of $c_1(L_k) = F(a,b)$. Comparing the second
Chern class gives, \begin{align} \label{5th}ab + bc + ca = [2](m)\cdot[2](n).
   \end{align} Comparing the third Chern class gives \begin{align}
     \label{5th3}abc = 0. 
  \end{align}  As the fourth Chern class of the two vector bundles is zero, it
does not give any relation.

We will now prove a higher dimensional version of lemma \ref{rel}
to show that relations \ref{5th1} and \ref{5th3} are implied by
relation \ref{5th}.

\bl Let $MU^*(BGL(1)^3) \cong MU^*[[r,s,t]]$. Let $u = [-1](r +_F s +_F
t)$. Let $d_1 = r + s + t + u$, $d_2 = rs + rt + ru + st + su + tu$,
$d_3 = rst + stu + rtu + rsu$, $d_4 = rstu$. Then there exist formal
power series $P_1(-,-)$ and $P_3(-,-)$ such that
\begin{align*}
  d_1 = P_1(d_2,d_4) \text{ and } d_3 = P_3(d_2,d_4). 
\end{align*}

\el

\pf

 Consider the map $r: GL(1)^3 \ra Sp(2)$ given by 

\begin{align}(u,v,w) \mapsto \begin{pmatrix}u & & & \\ & v & & \\ & & w & \\ & & &
  (uvw)^{-1}\end{pmatrix}\end{align}   

This induces the map $r^*: MU^*[[c_2,c_4]] \ra MU^*[[y_1,y_2,y_3]].$ Let the
standard representation of $Sp(2)$ be $\theta$. We know that
$MU^*(BSp(2) \cong MU^*[[c_2(\theta),c_4(\theta)]]$. Let $c_1(\theta) =
P_1(c_2(\theta),c_4(\theta))$ and $c_3(\theta) =
P_3(c_2(\theta),c_4(\theta))$.The pullback of $\theta$
under $r$ is $ r^*\theta = \phi_1 \oplus \phi_2 \oplus \phi_3 \oplus
\varphi$ where 
$\phi_i$ is the projection $GL(1)^3 \ra GL(1)$ and $\varphi$ is the dual
of $\phi_1 \otimes \phi_2 \otimes \phi_3$.  

We have $c_1(\phi_1) = r$, $c_1(\phi_2) = s$, $c_1(\phi_3) = t$ and
$c_1(\varphi) = u$. Then for each $i$, the expression $d_i$ is the
$i$th Chern of $r^*\theta$. Hence we have,   
\begin{align*}
  d_1 = P_1(d_2,d_4) \text{ and } d_3 = P_3(d_2,d_4). 
\end{align*}





\qed 

Hence we get
\begin{align*}
  a + b + c = P_1(ab + bc + ca,0) \text{ and }& abc = P_3(ab + bc +
  ca,0)\\
[2](m) + [2](n) = P_1([2](m)\cdot [2](n),0) \text{ and }& 0 =
P_3([2](m)\cdot [2](n),0).  
\end{align*}

This implies that relation \ref{5th} implies relations
\ref{5th1} and \ref{5th3}.

The sixth and the final relation in the cobordism of $BQ$ is obtained as
follows. Tensoring the two 
sides of $L_i + L_j + L_k + \mathds{1} = V^2$ by $V$ gives that \bd 4V
= V^3.\ed The relation we need is $c_2(4V) = c_2(V^3).$  

 Let $\eta$ denote the
standard representation $SL(2)$. We have $MU^*(BSL(2)) \cong MU^*[[e_2]]$. So let $c_2(4\eta) = G(e_2)$ and
$c_2(\eta^3) = H(e_2)$ where $G(-)$ and $H(-)$ are formal power
series. Pulling back $G(-)$ and $H(-)$ to $MU^*(BQ)$ via $Q \subset
SL(2)$, we have \begin{align}\label{6th} G(c_2(V)) = H(c_2(V)) 
\end{align} 

This gives us the six relations, namely, \ref{1st}, \ref{2nd},
\ref{3rd}, \ref{4th}, \ref{5th} and \ref{6th}. As all these relations
are symmetric in $m$ and $n$, we can 
write them in terms of $a$, $b$, $mn$ and $m+n$. Then using $m+n =
P(mn)$ we can write these relations in terms of $a$, $b$ and
$mn$. Taking their images in $MU^*(BQ)$, i.e. replacing $a$ by $x$,
$b$ by $y$ and $mn$ by $z$, we get the six equations in $MU^*(BQ)$. These
are precisely the six power series which generate the ideal $I$ in
$MU^*[[x,y,z]]$  so that $MU^*(BQ) \cong
MU^*[[x,y,z]]/I$. 

\bt $\o^*(BQ) \cong MU^*(BQ).$ \et

\pf We have already established $MU^*(BQ) \cong
MU^*[[x,y,z]]/I$ where $I$ is generated by the six relations mentioned above.
Observe that we can obtain these six relations in $\o^*[[x,y,z]]$ by
using exactly the same arguments as above because we already have
$MU^*(BG) \cong \o^*(BG)$ for $G= GL(n)$, $SL(n)$ and $Sp(n)$. And now
we go back to the argument in section \ref{lie}. Under the map
$\o^*(BQ) \ra 
MU^*(BQ)$ a relation $f(c_1(L_i), c_1(L_j), c_2(V)) = 0$ in $\o^*(BQ)$
is mapped to a relation $f(c_1(L_i), c_1(L_j), c_2(V))  = 0$ in
$MU^*(BQ)$. As there are no relations in $MU^*(BQ)$ other than the six
mentioned above, the same is true in $\o^*(BQ)$, hence \begin{align}
  \o^*(BQ) \cong MU^*(BQ) \end{align} 
\qed
  
\textbf{Remark :} It is complicated to obtain the
coefficients of these six power series. But we can at least see what
are the images of 
these relations in $CH^*(BQ) \cong H^*(BQ)$. The images of  $x$, $y$
and $z$ in $H^*(BQ)$ are $x'$, $y'$ and $z'$ respectively. The images
of relations \ref{1st} and \ref{2nd} are $x'^2 = y'^2 = 0$. The images
of relations \ref{3rd} and \ref{4th} are $2x' = 2y' = 0$. The
image of relation \ref{5th} is $x'y' = 4z'$. And finally the image
of relation \ref{6th} is $8z' = 0$. 

The Chow ring of the classifying space of an iterated
wreath product  $\Z/p \wr \Z/p \wr \cdots \wr \Z/p$ has been calculated
in \cite{tot99}. Let $D = D_8$ denote the dihedral group of order
$8$. Observe that $D \cong \Z/2 \wr \Z/2$.  Hence we have \bd
CH^*(BD) \cong \Z[x_D,y_D,z_D]/(2x_D = 2y_D = 4z_D =0; x_Dy_D =
2z_D)\ed where $x_D,y_D \in CH^1(BD)$ and $z_D \in CH^2(BD)$. Hence by 
lemma \ref{gen}, we have the structure of $\o^*(BD)$ as an $\o^*$-module.




\def\cprime{$'$} \def\cprime{$'$} \def\cprime{$'$}

\end{document}